\documentstyle[12pt, twoside]{article}
\input amssym.def
\input amssym.tex

\newenvironment{demo}[1]%
{\vskip-\lastskip\medskip
  \noindent
  {\em #1.}\enspace
  }%
{\qed\par\medskip
  }

\newcommand{\qed}{
  \strut\hfill
  \mbox{$\Box$}
  }

\newcommand{\C}{ {\Bbb C} }
\newcommand{\Ind}{ \mbox{Ind} }
\newcommand{\Res}{ \mbox{Res} }
\newcommand{\Z}{ {\Bbb Z} }

\newcommand{\FG}{ {\cal F}_G }

\newtheorem{theorem}{Theorem}

\newtheorem{lemma}{Lemma}

\newtheorem{remark}{Remark}

\newtheorem{definition}{Definition}

\newtheorem{proposition}{Proposition}

\newtheorem{corollary}{Corollary}

\begin{document}
\title{Equivariant K-theory, wreath products,
and Heisenberg algebra
}
\author{
  Weiqiang Wang\thanks{1991 {\em Mathematics Subject Classification}.
 Primary 19L47, 17B65.
        }
}
\date{}
\maketitle

\begin{abstract}{Given a finite group $G$ and a $G$-space $X$,
we show that a direct sum
${\cal F}_G (X) = \bigoplus_{ n \geq 0}K_{G_n} (X^n) \bigotimes \C$
admits a natural graded Hopf algebra and $\lambda$-ring structure,
where $G_n$ denotes the wreath product $G \sim S_n$.
$\FG (X)$ is shown to be isomorphic to a certain
supersymmetric product in terms
of $K_G(X)\bigotimes \C$ as a graded algebra.
We further prove that $\FG (X)$ is isomorphic to the Fock space of
an infinite dimensional Heisenberg (super)algebra.
As one of several applications, we compute the orbifold Euler characteristic
$ e( X^n, G_n)$.}
\end{abstract}

\setcounter{section}{-1}
\section{Introduction}
 Given a finite group $G$ and a locally
compact, Hausdorff, paracompact $G$-space $X$,
the $n$-th direct product $X^n$ admits a natural action
of the wreath product $G_n = G \sim S_n$ which is a semi-direct product
of the $n$-th direct product $G^n$ of $G$ and the
symmetric group $S_n$. The main goal of the present paper is to
study the equivariant topological K-theory $K_{G_n} (X^n)$,
for all $n$ together, and discuss several applications which are
of independent interest.

We will first show that a direct sum
$$\FG (X) = \bigoplus_{ n \geq 0}K_{G_n} (X^n) \bigotimes \C
$$
carries several wonderful structures.
More explicitly,
we show that $\FG (X)$ admits a natural Hopf algebra structure
with a certain induction functor as multiplication and a certain
restriction functor as comultiplication (cf. Theorem~\ref{th_hopf}).
When $X$ is a point, $K_{G_n} (X^n)$ is the Grothendieck
ring $R(G_n)$, and we recover the standard Hopf algebra
structure of $\bigoplus_{n \geq 0} R(G_n)$
(cf. e.g. \cite{M1, M2, Z}). A key lemma used here is a straightforward
generalization to equivariant K-theory
of a statement in the representation theory of
finite groups concerning the restriction of an induced
representation to a subgroup.

We show that $\FG (X)$ is a free $\lambda$-ring generated by
$K_G(X)\bigotimes \C$ (cf. Proposition~\ref{prop_ring}).
We write down explicitly  the Adams operations $\varphi^n$'s
in $\FG (X)$. Incidentally we also obtain an equivalent way
of defining the Adams operations in
$K_G(X) \bigotimes \C$ (not over $\Z$) by means of the wreath
products, generalizing a definition by Atiyah \cite{A} in terms
of the symmetric group
in the ordinary (i.e. non-equivariant) K-theory setting.
When $X$ is a point we recover the $\lambda$-ring structure
of $\bigoplus_{n \geq 0} R(G_n)$ (cf. \cite{M1}).

As a graded algebra $\FG (X)$
has a simple description as a certain supersymmetric algebra in terms
of $K_G(X)\bigotimes \C$ (cf. Theorem~\ref{th_main}).
The proof uses a theorem in \cite{AS} and the structures
of the centralizer group of an element in $G_n$ and of the
fixed point set of the action of $a \in G_n$ on $X^n$ which we
work out in Sect.~\ref{sec_wreath}.
In particular, this description indicates that  $\FG (X)$ has the size of
a Fock space of a certain infinite-dimensional Heisenberg superalgebra
which we will construct in terms of natural
additive maps in K-theory (cf. Theorem~\ref{th_heisenberg}).

Our results above generalize Segal's work \cite{S2},
and our proofs are direct generalizations of those in
\cite{S2} (also see \cite{Z, M1}).
What Segal studied in \cite{S2}, partly motivated by
remarks in Grojnowski \cite{Gr}, is the space
$\bigoplus_{ n \geq 0}K_{S_n} (X^n) \bigotimes \C$
for compact $X$ which corresponds to our special case when $G$ is trivial
and then $G_n$ is the symmetric group $S_n$.
The paper \cite{Gr} was in turn
motivated by a physical paper of Vafa and Witten \cite{VW}.
Our present work grew out of an attempt
to understand Segal's outlines \cite{S2} and was also stimulated by
Nakajima's lecture notes on Hilbert schemes \cite{N}
(also cf. \cite{N1, Gr}).
Our first main observation in this paper is
that there is a natural way to add the group $G$ into
Segal's scheme and this allows several different applications
as discussed below. These applications are of independent
interest on their own. We expect that there is also a natural
way to incorporate $G$ into the remaining part of \cite{S2}.

Our addition of $G$ has already highly non-trivial consequences
even in the case when $X$ is a point.
By tensoring ${\cal F}_G (pt)$ with the group algebra
of the Grothendieck ring $R(G)$, we obtain the underlying vector
space for the vertex algebra
associated to the lattice $R(G)$ \cite{B, FLM}.
When $G$ is a finite subgroup of $SL_2(\Bbb C)$, this
will lead to a group theoretic construction of
the Frenkel-Kac-Segal vertex
representation of an affine Kac-Moody Lie algebra, which can be
viewed as a new form of McKay correspondence \cite{Mc}.
Detail along these lines will be developed in a forthcoming paper.

An interesting case of our study is that $X $ is the complex plane $\C^2$
acted upon by a finite subgroup $G$ of $SL_2(\Bbb C)$.
Let $\widetilde{\C^2 /G}$ denote
the minimal resolution of singularities of $\C^2 /G$
(cf. e.g. \cite{N}). Via the McKay correspondence \cite{Mc}, we show that
either $K_{G_n} ((\C^2)^n) \otimes \C$
or $K_{S_n} ((\widetilde{\C^2 /G})^n) \otimes \C$
has the same dimension as the homology group of
the Hilbert scheme of $n$ points on $\widetilde{\C^2 /G}$
(cf. \cite{G}).
This fact has a straightforward generalization,
cf. Remark~\ref{rem_dim}.
Our message here is that the wreath product plays an important
role in the study of the Hilbert scheme of
$n$ points on $\widetilde{\C^2 /G}$ in exactly the way
a symmetric group $S_n$ does for the Hilbert scheme of
$n$ points on ${\C^2 }$ (cf. \cite{N, BG}), which
is in turn a special case of the former when $G$ is trivial.
We will discuss these in more detail in another occasion.

For a smooth manifold $X$ acted upon by $G$, Dixon, Harvey,
Vafa and Witten \cite{DHVW}
introduced a notion of orbifold Euler characteristic $e(X, G)$
in their study of string theory of orbifolds. We show that
the orbifold Euler characteristic $ e(X^n, G_n)$
is uniquely determined by $e(X, G)$ and $n$.
In terms of a generating function,
our formula reads (see Theorem~\ref{th_orbi}):
\begin{eqnarray}  \label{eq_formu}
  \sum_{ n \geq 1} e(X^n, G_n) q^n
   =  \prod_{ r =1}^{\infty} (1 -q^r)^{- e(X, G)}.
\end{eqnarray}
By putting $G =1$ and thus $e(X, G) = e(X)$,
we recover a formula of Hirzebruch-H\"ofer \cite{HH}.
By using Eq.~(\ref{eq_formu}) and G\"ottsche's
formula \cite{G}, we show that $X^n /G_n$
admits a resolution of singularities whose Euler characteristic coincides
with the orbifold Euler characteristic $e(X^n, G_n)$
assuming that $X$ is a smooth quasi-projective surface and $X/G$ has
a resolution of singularities whose Euler characteristic is $e (X, G)$.

In this paper the language of equivariant K-theory is used.
We should also mention the very relevant construction of
Heisenberg superalgebra on a direct sum over $n$
of the homology group $H(X^{[n]})$ of Hilbert
scheme of $n$ points on a smooth quasi-projective surface $X$,
due to Nakajima \cite{N1} and Grojnowski \cite{Gr} independently. However
the constructions and computations in terms of K-theory are
simpler and work for more general spaces.  Bezrukavnikov and Ginzburg
\cite{BG} have proposed a way to obtain a direct isomorphism
from $K_{S_n} (X^n) \bigotimes \Bbb C$ to $H(X^{[n]})$ for an algebraic
surface $X$. Independently
M.~de Cataldo and L.~Migliorini has recently established this isomorphism
for complex surfaces \cite{CM}.

The plan of this paper is as follows. In Sect.~\ref{sec_wreath}
we give a presentation of the wreath product $G_n$ and
study its action on $X^n$. In Sect.~\ref{sec_hopf} we
construct a Hopf algebra structure on $\FG (X)$.
In Sect.~\ref{sec_main} we give a description of $\FG (X)$
as a graded algebra.
In Sect.~\ref{sec_ring}
we give a $\lambda$-ring structure on $\FG (X)$.
In Sect.~\ref{sec_heis} we
construct the Heisenberg superalgebra which
acts on $\FG (X)$ irreducibly. In Sect.~\ref{sec_orbi} we calculate
the orbifold Euler characteristic $e(X^n, G_n)$
and study in detail the special case when
$X$ is the complex plane acted upon by a finite subgroup of $SL_2 ( \C )$.
We have included some detail which are probably
trivial to experts hoping that this may benefit readers
with different backgrounds.
\section{The wreath product and its action on $X^n$}
   \label{sec_wreath}
 Let $G$ be a finite group.
We denote by $G_*$ the set of conjugacy classes of $G$ and
$R(G)$ the Grothendieck ring of $G$. $R(G) \bigotimes_{\Z} \C$
can be identified with the ring of class functions $C(G)$ on $G$
by taking the character of a representation.
Denote by $\zeta_c$ the order of the centralizer
of an element lying in the conjugacy class
$c$ in $G$. We define an inner product on $C(G)$ as usual:
\begin{eqnarray}   \label{eq_inner}
( \chi |\psi ) = \frac1{|G|} \sum_{g \in G} \chi (g) \overline{\psi}(g),
\quad \chi, \psi \in C(G).
\end{eqnarray}

Let $G^n = G \times \ldots \times G$ be the direct product of
$n$ copies of $G$. Denote by $|G|$ the order of $G$, and by
$[g]$ the conjugacy class of $g \in G$.
The symmetric group $S_n$ acts on $G^n$
by permuting the $n$ factors:
$ s (g_1, \ldots, g_n) = (g_{s^{-1}(1)} , \ldots, g_{s^{-1}(n)} )$.
The {\em wreath product} $G_n = G \sim S_n$ is defined to be
the semidirect product of $G^n$ and $S_n$, namely the multiplication on
$G_n$ is given by
$(g, s)(h, t) = (g. s(h), st)$, where $g, h \in G^n, s, t \in
S_n$. Note that $G^n$ is a normal subgroup of
$ G_n$ by identifying $g \in G^n$ with
$(g,1) \in G_n$.

Given $a = (g, s) \in G_n$ where $g = (g_1, \ldots, g_n)$,
we write $s \in S_n$ as a product
of disjoint cycles: if  $z= (i_1, \ldots, i_r)$ is one of them,
the {\em cycle-product} $g_{i_r} g_{i_{r-1}} \ldots g_{i_1} $
of $a$ corresponding to the cycle $z$
is determined by $g$ and $z$ up to conjugacy.
For each $c \in G_*$ and each integer $r \geq 1$, let
$m_r (c)$ be the number of $r$-cycles in $s$ whose cycle-product
lies in $c$. Denote by $\rho (c)$ the partition having
$m_r (c)$ parts equal to $r$ ($r \geq 1$) and
denote by $\rho = ( \rho (c) )_{c \in G_*}$
the corresponding partition-valued
function on $G_*$. Note that
$|| \rho || : = \sum_{c \in G_*} |\rho (c)|
= \sum_{c \in G_*, r \geq 1} r m_r (c) = n$,
where $| \rho (c) |$ is the size of the partition $\rho (c)$.
Thus we have defined a map from $G_n$ to ${\cal P}_n (G_*)$,
the set of partition-valued
function $\rho = ( \rho (c) )_{c \in G_*}$ on $G_*$
such that $|| \rho || = n$ . The function $\rho$ is
called the {\em type} of $a = (g, s) \in G_n$.
Denote ${\cal P} (G_*) = \sum_{n \geq 0} {\cal P}_n (G_*)$.

Given a partition
$\lambda$ with $m_r$ $r$-cycles ($r \geq 1$), define
$z_{\lambda} = \prod_{r \geq 1} r^{m_r} m_r !$.
This is the order of the centralizer in $S_n$ of
an element of cycle-type $\lambda$.
We shall denote by $l(\lambda) = \sum_{ r \geq 1} m_r$
the length of $\lambda$.

Given a partition-valued function $\rho \in {\cal P}(G_*)$, we define
$l (\rho) = \sum_{ c \in G_*} l ( \rho (c))$ and
$$ Z_{\rho} = \prod_{c \in G_*} z_{ \rho (c)} \zeta_c^{l (\rho(c))}.$$
Denote by $\sigma_n (c)$ the class function of $G_n$ which takes
value $n \zeta_c$ at an $n$-cycle whose cycle-product lies in $c \in G_*$
and $0$ otherwise. For $\rho = \{m_r (c) \}_{c,r} \in  {\cal P}(G_*)$,
we define
$$\sigma^{\rho} = \prod_{c \in G_*, r \geq 1} \sigma_r(c)^{m_r(c)}.
$$
We regard $\sigma^{\rho}$ as the class function on
$G_n$ which takes value $Z_{\rho}$ at elements of
type $\rho$ (where $n = || \rho ||$) and $0$ elsewhere.

We formulate some well-known facts below (cf. e.g. \cite{M2}) which
will be needed later.

\begin{proposition}  \label{prop_type}
 Two elements in $G_n$ are conjugate to each
other if and only if they have the same type.
The order of the centralizer in $G_n$ of
an element of type $\rho$ is $Z_{\rho}$.
\end{proposition}

We want to calculate the centralizer $Z_{G_n} (a)$ of $a \in G_n$.
First we consider the typical case that $a$ has one $n$-cycle.
As the centralizers of conjugate elements are conjugate subgroups,
we may assume that $a$ is of the form
$a = ( (g, 1, \ldots, 1), \tau)$
by Proposition~\ref{prop_type}, where $ \tau = (1 2 \ldots n)$.
Denote by $Z_G^{\Delta}(g)$, or $Z_G^{\Delta_n}(g)$ when it is necessary
to specify $n$, the following diagonal subgroup of $G^n$
(and thus a subgroup of $G_n$):
\begin{eqnarray*}
 Z_G^{\Delta}(g) = \left\{ ( (h, \ldots, h), 1) \in G^n \mid h \in Z_G(g)
                 \right\}.
\end{eqnarray*}

The following lemma follows from a direct computation.
\begin{lemma} \label{lem_perm}
  The centralizer $Z_{G_n} (a)$ of $a$ in $G_n$ is equal to the product
 $Z_G^{\Delta}(g) \cdot \langle a \rangle$, where $\langle a \rangle$
 is the cyclic subgroup of $G_n$ generated by $a$.
 Moreover, $a^n \in Z_G^{\Delta}(g) $ and $|Z_{G_n} (a)| = n |Z_G(g)|.$
\end{lemma}

Take a generic element $a  = ( g, s) \in G_n$
of type $\rho = ( \rho (c) )_{c \in G_*}$, where
$\rho (c) $ has $m_r (c)$ $r$-cycles ($r \geq 1$).
By Proposition~\ref{prop_type},
we may assume (by taking a conjugation if necessary)
that the $m_r (c)$ $r$-cycles are of the form
$$
 g_{ur}(c) = ( (g, 1, \ldots,1), (i_{u1}, \ldots, i_{ur}) ),
  1 \leq u \leq m_r (c), g \in c.
$$
Denote $ g_r (c) = ( (g, 1, \ldots,1), (12 \ldots r ) ).$
Throughout the paper, $\prod_{c,r}$ is understood as the
product $\prod_{c \in G_*, r \geq 1}$.

\begin{lemma} \label{lem_cent}
   The centralizer $ Z_{G_n} (a)$ of $a \in G_n$ is isomorphic to
 a direct product of the wreath products
 \begin{eqnarray}  \label{eq_centra}
  \prod_{c,r} \left(
                 Z_{G_r} ( g_r (c) )
     \sim S_{m_r (c)} \right).
 \end{eqnarray}
 Furthermore $Z_{G_r} ( g_r (c) )$ is isomorphic to
 $Z^{\Delta_r}_G (g ) \cdot \langle g_{r} (c) \rangle$.
\end{lemma}

\begin{demo}{Proof}
  It follows from the first  part of Lemma~\ref{lem_perm}
 that the centralizer $ Z_{G_n} (a)$ should contain
 a certain subgroup naturally isomorphic to (\ref{eq_centra}).
 By the second part of Lemma~\ref{lem_perm} we can count
 that the order of (\ref{eq_centra}) is equal to $Z_{\rho}$.
 The lemma now follows by
 comparing with the order of $ Z_{G_n} (a)$ given
 in Proposition~\ref{prop_type}.
\end{demo}

We will use $\star$ to denote the multiplication
in $C(G_n)$ which corresponds to the tensor product in $R(G_n)$.
We denote by $\underline{n}$
the trivial representation of $G_n$, and $\underline{1^n}$
the sign representation of $G_n$ in which $G^n$ acts trivially
and $S_n$ acts by $\pm 1$ depending on a permutation is
even or odd.
By abuse of notations, we also use the same symbols to
denote the corresponding characters as well.
The following lemma follows easily from the definitions.

\begin{lemma}  \label{lem_orth}
 \begin{enumerate}
  \item[1)] Given $\rho, \widetilde{\rho} \in {\cal P}_n (G_*)$,
  $
  \sigma^{\rho} \star\sigma^{\widetilde{\rho} }
   = \delta_{\rho, \widetilde{\rho} } Z_{\rho} \sigma^{\rho}.
  $
 In particular,
 \begin{eqnarray}
 \underline{n} &= & \sum_{|| \rho || = n} Z_{\rho}^{ -1} \sigma^{\rho},
         \label{eq_triv}     \\
 \underline{1^n} &= &
     \sum_{|| \rho || = n} (-1)^{ n -l(\rho)}Z_{\rho}^{ -1} \sigma^{\rho}.
   \label{eq_sign}
 \end{eqnarray}

  \item[2)] $\left(
                 \sigma^{\rho} | \sigma^{\widetilde{\rho} }
             \right)
     = \delta_{\rho, \widetilde{\rho} } Z_{\rho}. $ In other
 words, $ \sigma^{\rho}$ takes value $Z_{\rho}$
 at the elements in $G_n$
 of type $\rho$ and $0$ elsewhere.
 \end{enumerate}
\end{lemma}

For a $G$-space $X$, we define an action of
$G_n$ on $X^n$ as follows: given $ a = ( (g_1, \ldots, g_n), s)$,
we let
\begin{eqnarray} \label{eq_action}
 a . (x_1, \ldots, x_n)
  = (g_1 x_{s^{-1} (1)}, \ldots, g_n x_{s^{-1} (n)})
\end{eqnarray}
where $x_1, \ldots, x_n \in X$. Next we want to determine the
fixed point set $( X^n )^a$ for $a \in G_n$.
Let us first calculate in the typical case
$a = ( (g, 1, \ldots, 1), \tau) \in G_n$.
Note that the centralizer group $Z_G(g)$ preserves
the $g$-fixed point set $X^g$.

\begin{lemma} \label{lem_elem}
  The fixed point set is
 \begin{eqnarray*}
  ( X^n )^a = \left\{ (x, \ldots, x) \in X^n\mid x= g x
              \right\}
 \end{eqnarray*}
 which can be naturally identified with $X^g$.
 The action of $Z_{G_n} (a)$ on $(X^n)^a$ can be identified canonically with
 that of $Z_G (g)$ on $X^g$ together with the trivial action of
 the cyclic group $\langle a \rangle$ (cf. Lemma~\ref{lem_perm}).
 Thus
 $$ (X^n)^a / Z_{G_n} (a) \approx X^g / Z_G (g).
 $$
\end{lemma}

\begin{demo}{Proof}
   Let $(x_1, \ldots, x_n) $ be in the fixed point
 set $(X^n)^a$. By Eq.~(\ref{eq_action}) we have
 \begin{eqnarray*}
  (x_1, x_2, x_3, \ldots, x_n)
  &= & a. (x_1, x_2, x_3, \ldots, x_n) \\
  &= & (g x_n, x_1, x_2, \ldots, x_{n-1}).
 \end{eqnarray*}
 So all $x_i (i =1, \ldots, n)$ are equal to, say $x$,
 and $gx = x$. The remaining statements follow from
 Lemma~\ref{lem_perm}.
\end{demo}

All $Z_G (g)$ are conjugate and all
$X^g$ are homeomorphic to each other for different
representatives $g$ in a fixed conjugacy class $c \in G_*$.
Also the orbit space $X^g /Z_G (g)$ can be identified
with each other by conjugation for different
representatives of $g$ in $c \in G_*$.
We make a convention to denote $Z_G (g)$ (resp. $X^g$, $X^g /Z_G (g)$)
by $Z_G (c)$ (resp. $X^c$, $X^c /Z_G (c)$) by abuse of notations
when the choice of a representative $g$ in $c$ is immaterial.

\begin{lemma} \label{lem_fix}
 Retain the notations in Lemma \ref{lem_cent}.
 The fixed point set $( X^n )^a$ can be naturally identified
 with $ \prod_{c,r}  (X^c )^{m_r (c)} $.
 Furthermore the orbit space
 $ ( X^n )^a / Z_{G_n}(a)$ can be naturally
 identified with
 $$ \prod_{c,r} S^{m_r (c)} \left( X^{c} / Z_G (c) \right) $$
 where $S^{m}(\cdot)$ denotes the $m$-th symmetric product.
 \end{lemma}

\begin{demo}{Proof}
 The first part easily follows from Lemma~\ref{lem_elem}.
 By Lemma \ref{lem_cent} and Lemma \ref{lem_elem},
 the action of $Z_{G_n}(a)$ on $( X^n )^a$ can be naturally
 identified with that of
 $$
  \prod_{c,r} \left(
               \left(
                   Z^{\Delta_r}_G (g ) \cdot \langle g_{r} (c) \rangle
               \right)
                    \sim S_{m_r (c)}
              \right)
 $$
 on $ \prod_{c,r}  (X^c )^{m_r (c)} $
 where $S_{m_r (c)}$ acts by permutation and
 $\langle g_{r} (c) \rangle$ acts on $X^{c}$ trivially.
 Thus the second part of the lemma follows.
\end{demo}
\section{The Hopf algebra structure on $\FG (X)$}
   \label{sec_hopf}
Given a compact Hausdorff $G$-space $X$, we recall \cite{A1, S1} that
$K^0_G(X)$ is the Grothendieck group of
$G$-vector bundles over $X$. One can define $K_G^1 (X)$
in terms of the $K^0$ functor and a certain suspension operation,
and one puts
$$ K_G(X) =  K^0_G(X) \bigoplus K^1_G(X).
$$
The tensor product of vector bundles gives rise to
a multiplication on $ K_G(X)$ which is super (i.e. $\Z_2$-graded)
commutative. In this paper we will be only concerned about
the free part $K(X) \otimes \C$, which will be
denoted by $\underline{K}_G(X)$ subsequently.
We denote by $\dim {K}_G^i (X) (i = 0, 1)$
the dimension of ${K}_G^i (X) \otimes \C$.

If $X$ is locally compact, Hausdorff and
paracompact $G$-space, take the one-point
compactification $X^+$ with the extra point $\infty$ fixed by $G$.
Then we define $K_G^0 (X)$ to be the kernel of the map
$$ K^0_G(X^+) \longrightarrow K^0_G( \{\infty\} )
$$
induced by the inclusion map $\{\infty\} \hookrightarrow X^+$.
It is clear that this definition is equivalent to the
earlier one when $X$ is compact. We also define
$K^1_G(X) = K^1_G(X^+).$

Note that $K_{G} (pt)$ is isomorphic to the Grothendieck
ring $R(G)$ and $\underline{K}_{G} (pt)$ is
isomorphic to the ring $C(G)$ of class functions on $G$.
The bilinear map $\star$ induced from the tensor product
$$K_{G} (pt) \bigotimes K_{G} (X) \longrightarrow K_{G} (X) $$
gives rise to a natural $\underline{K}_{G} (pt)$-module structure
on $\underline{K}_{G} (X)$. Thus $\underline{K}_{G} (X)$
naturally decomposes into a direct sum over the set of
 conjugacy classes $G_*$ of $G$.
The following theorem \cite{AS} (also cf. \cite{BC} )
gives a description of each direct summand.

\begin{theorem} \label{th_tech}
  There is a natural $\Z_2$-graded isomorphism
  \begin{eqnarray*}
   \phi : \underline{K}_{G} (X)  \longrightarrow
    \bigoplus_{[g]} \underline{K} \left( X^{g} /Z_{G} (g)
                                  \right).
  \end{eqnarray*}
\end{theorem}

Given $c \in G_*$, we denote by $\sigma_c$ the
class function which takes value
$\zeta_c$ at an element of $G$ lying in the conjugacy
class $c$ and $0$ otherwise. Then an element
in $\underline{K}_{G} (X) $ can be written as of the form
$\sum_{c \in G_*} \xi_{c} \sigma_c$, where
$\xi_c \in \underline{K} \left( X^{g} /Z_G (g) \right).$
More explicitly the isomorphism $\phi$ is defined as follows
when $X$ is compact:
if $E$ is a $G$-vector bundle over $X$ its restriction to
$X^g$ is acted by $g$ with base points fixed.
Thus $E|_{X^g }$ splits into a direct sum of
subbundles $E_{\mu}$ consisting of eigenspaces of
$g$ fiberwise for each eigenvalue $\mu$ of $g$.
$Z_G (g)$ acts on $X^g$ and one may check that
$\sum_{\mu} \mu E_{\mu}$ indeed lies in the $Z_G (g)$-invariant
part of $\underline{K}^0 (X^g) $.
Put
$$
 \phi_{g} (E) =  \sum_{\mu} \mu E_{\mu}
 \in \underline{K}^0 \left( X^{g} /Z_G (g) \right)
  =  \underline{K}^0  \left( X^{g} \right)^{Z_G (g)}.
$$
The isomorphism $\phi$ on the $K^0$ part is given by
$ \phi = \bigoplus_{[g] \in G_*} \phi_{g}$.
Then one easily extends the isomorphism $\phi$ to $K^1$
as $K^1_G (X)$ can be identified with the kernel of the
map from $K^0_G (X \times S^1)$ to $K^0_G (X)$
given by the inclusion of a point in $S^1$.
When $X$ is a point the isomorphism $\phi$
becomes the map from a representation of $G$ to its character.
By some standard arguments using compact pairs \cite{A1}, the isomorphism
remains valid for a locally compact, Hausdorff and paracompact $G$-space.
The following lemma is well known.

\begin{lemma}
   Given a finite group $G$,
 a subgroup $H$ of $G$, and a $G$-space $X$. There is
 a natural {\em induction functor}
 $\Ind = \Ind_H^G : K_H (X) \longrightarrow K_G (X)$.
 In particular when $X$ is a point, the functor $\Ind_H^G$
 reduces to the familiar induction functor of representations.
\end{lemma}

\begin{demo}{Proof}
  Note that there is a $G$-equivariant isomorphism
 \begin{eqnarray}   \label{eq_equiva}
  G \times_H X \stackrel{\mbox{Iso}}{\longrightarrow} G/H \times X
 \end{eqnarray}
 by sending $(g, x) \in G \times_H X$ to $(gH, gx)$. We remark
 that although both sides of (\ref{eq_equiva}) remain well-defined
 for a $H$-spcase $X$ without a $G$-action, the map $\mbox{Iso}$
 makes sense only for a $G$-space $X$.
 Denote by $p: G\times_H X \longrightarrow X$ the
 composition of the projection $G/H \times X$ to $X$ with
 the isomorphism (\ref{eq_equiva}).
 As is well known \cite{S1}, one has a natural isomorphism
 $$K_H (X) \longrightarrow K_G(G \times_H X)$$
 by sending an $H$-equivariant vector bundle $V$ on $X$
 to the $G$-equivariant vector bundle $G \times_H V$.
 The composition $p \circ \pi$
 of the projection $p: G\times_H X \longrightarrow X$
 with the bundle map $\pi : G\times_H V \longrightarrow G\times_H X$
 sends $(g, v)$ to $(g, \pi (v) )$. One easily check that this gives rise to
 a well-defined $G$-equivariant vector bundle on $X$, which induces
 the induction functor
 $\Ind_H^G : K_H (X) \longrightarrow K_G (X)$.
\end{demo}

 We denote by $\Res^G_H$ (or $\Res_H$,
 or even $\Res$, if there is no ambiguity)
the {\em restriction functor} from $K_G(X)$ to $K_H(X)$
by regarding a $G$-equivariant vector bundle as
an $H$-equivariant vector bundle. Denote
\begin{eqnarray*}
 \FG (X) = \bigoplus_{n \geq 0} \underline{K}_{G_n} (X^n),
 \quad
 {\cal F}_G^q (X) = {\bigoplus}_{n \geq 0}
  q^n \underline{K}_{G_n} (X^n)
\end{eqnarray*}
where $q$ is a formal variable counting the graded structure of $ \FG(X)$.
We introduce the notion of $q$-dimension:
$$\dim_q \FG (X) = \sum_{n \geq 0} q^n \dim K_{G_n} (X^n).
$$

Define a multiplication $\cdot$ on $ {\cal F}_G (X)$ by a composition
of the induction map and the K\"unneth isomorphism $k$:
\begin{eqnarray}  \label{eq_mult}
   \underline{K}_{G_n} (X^n) \bigotimes  \underline{K}_{G_m} (X^m)
      \stackrel{k }{\longrightarrow} \underline{K}
           _{G_n \times G_m} (X^{n +m} )
 \stackrel{\mbox{Ind}}{\longrightarrow}  \underline{K}_{G_{n+m}} (X^{n +m} ).
\end{eqnarray}
We denote by $ 1 $ the unit in $K_{G_0} (X^0) \cong \C $.

On the other hand we can define a comultiplication
$\Delta$ on $ {\cal F}_G (X)$, given by a composition of
the inverse of the K\"unneth isomorphism and the
restriction from $G_n$ to $ G_k  \times G_{n -k}$:
\begin{eqnarray*}
   \underline{K}_{G_n} (X^n)
 \longrightarrow \bigoplus_{m =0}^n
        \underline{K}_{G_{m } \times G_{n -m} } (X^n)
 \stackrel{k^{-1}}{ \longrightarrow} \bigoplus_{m =0}^n
    \underline{K }_{G_{m }} (X^m) \bigotimes
    \underline{K}_{G_{n -m}} (X^{n -m}).
\end{eqnarray*}
We define the counit $\epsilon : {\cal F}_G (X) \longrightarrow \C$
by sending $\underline{K}_{G_n} (X^n)$ $(n >0)$ to $0$
and $ 1 \in K_{G_0} (X^0) \approx \C $ to $1$.
The antipode can be also easily defined (see Remark~\ref{rem_hopf}).

\begin{theorem} \label{th_hopf}
  With various operations defined as above,
 $ \FG (X)$ is a graded Hopf algebra.
\end{theorem}

 To prove Theorem \ref{th_hopf}, we will need some preparation.
Given two subgroups $H$ and $L$ of a finite group $\Gamma$
and a $\Gamma$-space $Y$,
and let $V$ be an $H$-equivariant vector bundle on $Y$.
We denote the action of $H$ on $V$ by $\rho$.
Choose a set of representatives
$S$ of the double coset $H \backslash \Gamma /L$.
$H_s = sHs^{-1} \cap L$ is a subgroup of $L$ for $s \in S$.
We denote by $V_s$ the $H_s$-equivariant vector bundle on $Y$
which is the same as $V$ as a vector bundle and has
the conjugated action
\begin{eqnarray}
 \rho^s (x) = \rho (s^{-1} x s), \;\;x \in H_s. \label{eq_conj}
\end{eqnarray}

\begin{lemma}  \label{lem_ser}
 $ \Res_L \Ind_H^{\Gamma} V$ is isomorphic to the direct sum
 of the $L$-equivariant vector bundles $\Ind^L_{H_s} V_s$ for all
 $s \in H \backslash \Gamma/ L$.
\end{lemma}

One easily shows that one can extend $V$ in Lemma~\ref{lem_ser} to
the whole $K_H(Y)$. In the case $Y = pt$,
an $H$-equivariant vector bundle is just an $H$-module,
and the induction and restriction functors become the more familiar
ones in representation theory. In such a case, the above lemma is
standard (cf. e.g. \cite{Ser}, Proposition 2.2). In view of
our construction of the induction functor and the restriction
functor, the proof of Lemma~\ref{lem_ser} is essentially
the same as in the case $X = pt$ which we refer to \cite{Ser}
for a proof.

\begin{demo}{Proof of Theorem \ref{th_hopf}}
  We will show below that the comultiplication $\Delta$ is
 an algebra homomorphism. The other Hopf algebra axioms are
 easy to check.

   We apply Lemma~\ref{lem_ser}
 to the case $Y = X^N$, $L = G_m \times G_n$, $H = G_l \times G_r$,
 and $\Gamma = G_N$, where $l + r = m + n = N$. In this
 case the double coset $H \backslash \Gamma/L$ is isomorphic
 to $(S_l \times S_r) \backslash S_N / (S_{m} \times S_{ n})$ since
 $G_N = G^N \cdot S_N$ and $G_l \times G_r = G^N \cdot (S_l \times S_r)$.
 Furthermore $(S_l \times S_r) \backslash S_N / (S_{m} \times S_{ n})$
 is parametrized by the $2 \times 2$ matrices
\begin{eqnarray}   \label{matr}
 \left[ \begin{array}{cc}
  a_{11}      & a_{12}      \\
  a_{21}      & a_{22}
 \end{array} \right]
\end{eqnarray}
 satisfying
 \begin{eqnarray}
  a_{ij} \in\Z_+,      & i, j = 1, 2,           \nonumber \\
  a_{11} + a_{12} =m , & a_{21} + a_{22} = n,   \label{eq_condition} \\
  a_{11} + a_{21} =l , & a_{12} + a_{22} = r.   \nonumber
 \end{eqnarray}
 We denote by $\cal M$ the set of all the $2 \times 2$ matrices
 of the form (\ref{matr}) satisfying the conditions (\ref{eq_condition}).

 Given $E \in \underline{K}_{G_l} (X^l), F \in \underline{K}_{G_r} (X^r),$
we calculate by using Lemma~\ref{lem_ser} as follows:

 \begin{eqnarray}  \label{eq_hopf}
   \Res_{(m,n)} \Ind_{(l, r)}^N (E \boxtimes F ) &= &
  \bigoplus_{A \in \cal M}
  \Ind^{(m,n)}_A (1_{a_{11} } \otimes T^{(a_{12}, a_{21})} \otimes 1_{a_{22} })
  (\Res_{A'} (E \boxtimes F)) \nonumber   \\
   &= &
   \bigoplus_{A \in \cal M} \Ind^m_{(a_{11}, a_{12})} F_1 \boxtimes
   \Ind^n_{(a_{21}, a_{22})} F_2 .
 \end{eqnarray}
 Here the superscript or subscript ${(a,b)}$
 is a short-hand notations for $G_a \times G_b$. $1_a$
 stands for the identity operator from $\underline{K}_{G_a} (X^a)$
 to itself, and $\Ind^{a +b}_{(a,b)}$ for the induction functor from
 $\underline{K}_{ G_a \times G_b} (X^{a +b})$
 to $\underline{K}_{ G_{a +b} } (X^{a +b})$. $T^{(a,b)}$
 denotes the canonical functor from
 $\underline{K}_{ G_a} (X^{a})\otimes \underline{K}_{ G_b} (X^{b})$
 to  $\underline{K}_{ G_b} (X^{b}) \otimes \underline{K}_{ G_a} (X^{a})$
 by switching the factors with an appropriate sign coming from the
 $\Z_2$ grading of K-theory. Given $A \in \cal M$ of the form
 (\ref{matr}), the $A$ in the expressions $\Res_A$,
 $\Ind_A$ etc stands for
 $G_A \equiv G_{a_{11}} \times G_{a_{12}} \times G_{a_{21}} \times G_{a_{22}}$
 while ${A'}$ for
 $G_{A'} \equiv
   G_{a_{11}} \times G_{a_{21}} \times G_{a_{12}} \times G_{a_{22}}$.
 We wrote $(1 \otimes T^{(a_{12}, a_{21})} \otimes 1)
  (\Res_{A'} (E \boxtimes F))$
 as $F_1 \boxtimes F_2$  instead of a direct
 sum of the form $F_1 \boxtimes F_2$ in order to simplify notations,
 with $F_i$ $ (i =1,2)$ the corresponding elements in
 $K_{G_{a_{i1} } \times G_{a_{i2} } }(X^{a_{11} +a_{12}} ) $.

 Now it is straightforward to check
 that the statement that $\Delta$ is an algebra homomorphism
 is just a reformulation of the identity obtained by summing
 Eq.~(\ref{eq_hopf}) over all possible $(m, n)$ with $m +n = N$.
\end{demo}

\section{A description of $\FG (X)$ as a graded algebra}
\label{sec_main}

In this section, we give an explicit description of $\FG (X)$
as a graded algebra which in particular tells us the
dimension of $\underline{K}_{G_n}(X^n)$.
\begin{theorem} \label{th_main}
  As a $(\Z_+ \times \Z_2)$-graded algebra
$ {\cal F}_G^q (X)$ is isomorphic to the supersymmetric algebra
  $ {\cal S} \left( \bigoplus_{ n \geq 1} q^n \underline{K}_G(X)
      \right)$. In particular,
 $$
 \dim_q {\FG} (X) =
 \frac{ \prod_{r \geq 1} (1 + q^r)^{ \dim K^1_G (X)} }{ \prod_{r \geq 1}
    (1 - q^r)^{ \dim K^0_G (X)} }.
 $$
\end{theorem}
The supersymmetric algebra here is equal to the tensor product
of the symmetric algebra
$ S \left( \bigoplus_{ n \geq 1} q^n \underline{K}^0_G(X)
      \right)$ and the exterior algebra
$\Lambda \left( \bigoplus_{ n \geq 1} q^n \underline{K}^1_G(X) \right)$.

\begin{demo}{Proof}
 Take $a \in G_n$ of type $\rho = \{m_r(c)\}_{c,r}$
 as in Sect.~\ref{sec_wreath}.
 By Lemma~\ref{lem_cent} and Lemma~\ref{lem_fix} and
 the K\"unneth formula, we have
 \begin{eqnarray}
   K( (X^n)^a / Z_{G_n} (a) )
   & \approx & \bigotimes_{c \in G_*, r \geq 1}
    \left(
      \left( K(X^{c})^{ Z_G (c)}
      \right)^{\bigotimes m_r (c)}
    \right)^{S_{m_r (c)}  }  \nonumber\\
   & \approx & \bigotimes_{c \in G_*, r \geq 1}
   {\cal S}^{m_r (c)} ( K(X^{c} / Z_G (c)) ).   \label{eq_symm}
 \end{eqnarray}

 Thus if we take a summation of Eq. (\ref{eq_symm})
 over all conjugacy classes of $G_n$
 and all over $n \geq 0$, we obtain:
 \begin{eqnarray*}
 {\cal F}^q_G (X)
   & \approx & \bigoplus_{ n \geq 0}
     \bigoplus_{ \{m_r(c)\}_{c,r} \in {\cal P}_n ( G_*)}
        q^n \bigotimes_{c, r}
        {\cal S}^{m_r (c)} (\underline{K} (X^{c} / Z_G (c)) )
       \quad\quad\quad\mbox{by Theorem } \ref{th_tech}, \label{eq_key}     \\
   & = & \bigoplus_{\{m_r(c)\}_{c,r} \in {\cal P} ( G_*)}
        \bigotimes_{c, r} {\cal S}^{m_r (c)}
          \left(q^r \underline{K}(X^{c} / Z_G (c)) \right)    \nonumber     \\
   & = & \bigoplus_{\{m_r\}_{r} } \bigotimes_{ r \geq 1} {\cal S}^{m_r}
          \left( \bigoplus_{c \in G_*} q^r \underline{K} (X^{c} / Z_G (c))
          \right)
       \mbox{    by letting $m_r = \sum_{ c \in G_*} m_r (c)$},
                               \label{eq_sum}     \\
   & = & \bigoplus_{\{m_r\}_{r} } \bigotimes_{ r \geq 1}
             {\cal S}^{m_r} (q^r \underline{K}_G (X) )
    \; \quad\quad\quad\quad\quad\quad\quad\quad\quad
    \quad\quad\quad\mbox{by Theorem } \ref{th_tech}, \label{eq_again}       \\
   & = & {\cal S} \left( \bigoplus_{ r \geq 1} q^r \underline{K}_G(X) \right).
                  \nonumber
 \end{eqnarray*}

  The statement concerning $\dim_q {\FG} (X)$ is an immediate consequence.
\end{demo}

\begin{remark} \rm  \label{rem_spec}
 \begin{enumerate}
 \item[1)] Theorem \ref{th_main} in the case when
 $G$ is the trivial group (and so $G_n = S_n$)
 is due to Segal \cite{S2}. Our proof is adapted from his to the
 wreath product setting.

 \item[2)] If $G$ acts on $X$ freely, so does $G^n$ on $X^n$.
 Then we have the isomorphism $K (X/G) \approx K_G (X) $. Note that
 $ G^n$ is a normal subgroup of the wreath product $G_n$
 and the quotient $G_n / G^n$ is isomorphic to $S_n$.
 By Proposition 2.1 in Segal \cite{S1} we see that
 \begin{eqnarray} \label{eq_free}
  K_{G_n} (X^n) \approx
  K_{G_n / G^n} (X^n / G^n) = K_{S_n} ( (X/G)^n).
 \end{eqnarray}
 Therefore Theorem \ref{th_main} follows from the special case
 $G = 1$ of Theorem \ref{th_main} (applying to $X/G$) and Eq. (\ref{eq_free}).

 \item[3)] When $X$ is a point, $\FG (pt) = {\bigoplus}_{n \geq 0} R( {G_n})$,
  and $\sigma^{\rho}, \rho \in {\cal P}( G_*)$
 form a linear basis for $\FG (pt)$ (cf. \cite{M2}). In particular,
 $$\dim_q \FG (pt) =  \prod_{r \geq 1} (1 - q^r)^{ -|G_*|}.
 $$

  \item[4)]  One may reformulate Theorem~\ref{th_main} in terms of the
 de-localized equivariant cohomology \cite{BBM, BC} via the Chern character.
  \end{enumerate}
\end{remark}

\begin{remark} \label{rem_hopf} \rm
  The Hopf algebra defined in Sect.~\ref{sec_hopf}
 can be identified via the isomorphism in Theorem~\ref{th_main}
 with the standard one on the supersymmetric algebra
 $ {\cal S} \left( \bigoplus_{ n \geq 1} \underline{K}_G(X) [n]
      \right)$
 by showing the sets of primitive vectors correspond to each other.
 Here $\underline{K}_G(X) [n]$ denotes the
 $n$-th copy of $\underline{K}_G(X)$ (see Theorem~\ref{th_main}). The
 antipode of the former space can be transfered via the isomorphism
 from the latter one.
\end{remark}

\section{The $\lambda$-ring structure on $\FG (X)$}
   \label{sec_ring}
Let us denote by $c_n$ the conjugacy class in $G_n$ which
has the type of an $n$-cycle and whose cycle product lies in the
conjugacy class $c \in G_*$. We consider the following
diagram of K-theory
 maps:

 \begin{eqnarray*}
  \underline{K}_G (X)
    & \stackrel{\boxtimes n}{\longrightarrow} & \underline{K}_{G_n} (X^n)
     \stackrel{\phi_n}{\longrightarrow}
    \bigoplus_{[a] \in (G_n)_*}
      \underline{K} \left( (X^n)^a / Z_{G_n} (a) \right) \\
     & \stackrel{\stackrel{pr}{\rightleftarrows} }{\iota}&
    \bigoplus_{c \in G_*}
      \underline{K} \left( (X^n)^{c_n} / Z_{G_n} (c_n) \right)
    \stackrel{\vartheta}{\longrightarrow}
    \bigoplus_{c \in G_*} \underline{K} \left( X^{c} / Z_G (c) \right)   \\
    & \stackrel{\phi}{\longleftarrow} &  \underline{K}_G (X) .
 \end{eqnarray*}

Given a $G$-equivariant vector bundle $V$,
we define a $G_n$-action on the $n$-th outer
tensor product $V^{\boxtimes n}$ by letting
\begin{eqnarray}  \label{eq_act}
 ( (g_1, \ldots, g_n), s). v_1 \otimes \ldots \otimes v_n
 = g_1 v_{ s^{-1}(1)} \otimes \ldots \otimes  g_n v_{ s^{-1}(n)}
\end{eqnarray}
where $g_1, \ldots, g_n \in G, s \in S_n$.
Clearly $V^{\boxtimes n}$ endowed with such a $G_n$ action
is a $G_n$-equivariant vector bundle over $X^n$.
Sending $V$ to $V^{\boxtimes n}$ gives rise to
the K-theory map $\boxtimes n$. $\phi_n$ is the isomorphism
in Theorem~\ref{th_tech} when applying to the case $X^n$ with the
action of $G_n$. $pr$ is the projection to the direct sum
over the conjugacy classes of $G_n$ which are of the type of an
$n$-cycle while $\iota$ denotes the inclusion map. $\vartheta$
denotes the natural identification given by Lemma~\ref{lem_fix}.
Finally the last map $\phi$ is the isomorphism given
in Theorem~\ref{th_tech}.

We shall now define a $\lambda$-ring
structure on ${\cal F}_G (X)$. It suffices to define the
Adams operations $\phi^n$ on ${\cal F}_G (X)$. We also
introduce several other K-theory operations which will
be needed later.

\begin{definition}
  We define the following composition maps:
 \begin{eqnarray*}
  \psi^n &:= &
  n \phi^{-1} \circ \vartheta \circ pr \circ \phi_n \circ \boxtimes n:
    \underline{K}_G(X) \longrightarrow  \underline{K}_G(X),  \\
 \varphi^n &:= &
  n \phi_n^{-1} \circ \iota \circ {pr} \circ \phi_n \circ  \boxtimes n:
    \underline{K}_G(X) \longrightarrow \underline{K}_{G_n}(X^n), \\
 {ch}_n &:= &  \phi^{-1} \circ \vartheta \circ pr \circ \phi_n:
   \underline{K}_{G_n}(X^n) \longrightarrow \underline{K}_G(X),  \\
 {\omega}_n &:= & n \phi_n^{-1} \circ \iota \circ \vartheta^{-1} \circ \phi:
    \underline{K}_G(X) \longrightarrow \underline{K}_{G_n}(X^n).
 \end{eqnarray*}
\end{definition}

We list some properties of these K-theory maps whose proof is straightforward.

\begin{proposition} \label{prop_property}
 The following identities hold:
 \begin{eqnarray*}
  ch_n \circ \omega_n  = n \;Id,\quad
  \omega_n \circ \psi^n  = n \varphi^n,\quad
  ch_n \circ \varphi^n  = n \psi^n.
 \end{eqnarray*}
\end{proposition}

Recall that $\sigma_n (c)$ is the class function of $G_n$ which takes
value $n \zeta_c$ at an $n$-cycle whose cycle-product lies in $c \in G_*$
and $0$ otherwise.
\begin{lemma} \label{lem_iden}
 \begin{enumerate}
  \item[1)]  $\varphi^n (V) =
       \sum_{c \in G_*} \zeta_c^{-1} V^{ \boxtimes n} \star \sigma_n (c).$
  \item[2)]   Both $\psi^n$ and $\varphi^n$ are additive K-theory maps.
 \end{enumerate}
\end{lemma}

  Note that the order of the centralizer of an element
lying in the conjugacy class $c_n$  is equal to $n \zeta_c$.
The first part of the above lemma
now follows from the definition of $\varphi^n $ and
Lemma~\ref{lem_orth}. The second part can be proved in exactly
the same way as in the symmetric group case \cite{A}. We record
here only a useful formula in the proof:
 Given $V, W $ two $G$-equivariant vector bundles,
let $[V]$ denote the corresponding element of $V$ in $K_G(X)$.
 (In general we use $V$ itself to denote the corresponding
 element in $K_G(X)$ by abuse of notation). Then
 \begin{eqnarray}
  ([V] -[W])^{\boxtimes n} = \sum_{j =0}^n (-1)^j
  \Ind_{G_{n -j} \times G_j}^{G_n}
  [V^{\boxtimes (n -j)} \boxtimes W^{\boxtimes j} ] \in K_{G_n}(X^n).
 \end{eqnarray}
 Here $V^{\boxtimes (n -j)}$ endows the standard $G_{n -j}$-action
 given by substituting $n$ with $n -j$ in Eq. (\ref{eq_act}), and
 $G_j$ acts on $W^{\boxtimes j}$ by the tensor product of the
 standard $G_j$-action with the sign representation of $G_j$.

$\psi^n$'s are the Adams operations on $\underline{K}_G (X)$,
giving rise to the $\lambda$-ring structure on $\underline{K}_G (X)$.
Theorem \ref{th_main} ensures us that $\FG (X)$ as a $\lambda$-ring
is free and generated by $\underline{K}_G(X)$.
\begin{proposition}   \label{prop_ring}
 $\FG (X)$ is a free $\lambda$-ring generated
 by $K_G(X) \bigotimes \C$,
 with $\varphi^n$'s as the Adams operations.
\end{proposition}

\begin{remark} \rm
 If $X$ is a point, then $K_{G_n} (pt) = R(G_n)$ and
 $\FG (pt) = \bigoplus_{n \geq 0} R(G_n).$ Our result reduces to
 the fact that $\FG (pt)$ is a free $\lambda$-ring
 generated by $G_*$ \cite{M1}. In the case when $G =1$, the
 proposition was due to Segal \cite{S2}.
\end{remark}

Denote by $\widehat{\cal F}_G^q (X)$ the completion
of $ {\cal F}_G^q (X)$ which allows formal infinite sums.
Given $V \in \underline{K}_G(X)$,
we introduce $H(V, q), E(V, q) \in \widehat{\cal F}_G^q (X)$ as follows:
\begin{eqnarray*}
 H(V, q) & =& \bigoplus_{n \geq 0} q^n V^{\boxtimes n},   \\
 E(V, q) & =& \bigoplus_{n \geq 0}
             q^n V^{\boxtimes n} \star \underline{1^n }.
 \end{eqnarray*}
\begin{proposition}
  One can express
 $H(V, q)$ and $E(V, q)$ in terms of $\varphi^r (V)$ as follows:
\begin{eqnarray}
 H(V, q) & =& \exp \left( \sum_{r >0} \frac1r \varphi^r (V) q^r \right),
   \label{eq_vo}  \\
 E(V, -q) & =& \exp \left(- \sum_{r >0} \frac1r \varphi^r (V) q^r \right).
   \nonumber
 \end{eqnarray}
\end{proposition}

\begin{demo}{Proof}
 We shall prove Eq.~(\ref{eq_vo}) by using Eq.~(\ref{eq_triv}).
 The formula for $E(V, -q)$ can be similarly obtained by using
 Eq.~(\ref{eq_sign}).
 \begin{eqnarray*}
 H(V, q) & =& \bigoplus_{n \geq 0} q^n V^{\boxtimes n} \star
           \underline{n}   \\
   & =& \bigoplus_{n \geq 0} q^n V^{\boxtimes n} \star
     \left( \sum_{|| \rho || = n} Z_{\rho}^{ -1} \sigma^{\rho}
     \right)                    \\
   & =& \bigoplus_{n \geq 0} \sum_{|| \rho || = n}
      \left( Z_{\rho}^{ -1} q^n V^{\boxtimes n} \star \sigma^{\rho}
      \right)                   \\
   & \stackrel{(\star)}{=}& \prod_{c \in G_*, r \geq 1}\frac1{m_r (c)!}
      \left( \frac1r q^r \zeta_c^{ -1} V^{\boxtimes r} \star \sigma_r (c)
      \right)^{m_r (c)}                                    \\
   & =&  \exp \left( \sum_{ r \geq 1} \frac 1r q^r
                   \sum_{c \in G_*} \zeta_c^{ -1}
                     V^{\boxtimes r} \star \sigma_r (c)
              \right) \quad\quad   \mbox{  by Lemma }  \ref{lem_iden},    \\
   & =&  \exp \left(  \sum_{ r \geq 1} \frac 1r q^r \varphi^{r} (V)
              \right) .
 \end{eqnarray*}
 Here the equation ($\star$) is understood by means of the multiplication
 in $\FG (X)$ given by the composition (\ref{eq_mult}).
\end{demo}

\begin{corollary}
   The $\lambda$-operations $\lambda^n$ on the $\lambda$ ring
 $\FG (X)$ sends $V \in \underline{K}_G(X)$ to
 $V^{\boxtimes n} \star \underline{1^n }$.
\end{corollary}

Combining with the additivity of $\varphi^r$, we have the following.
\begin{corollary}
 The following equations hold for $V, W \in \underline{K}_G(X)$:
 \begin{eqnarray*}
   H( -V, q) & =& E(V, -q)  \\
  H( V \bigoplus W, q) & =& H( V, q)H( W, q).
 \end{eqnarray*}
\end{corollary}

\section{$\FG (X)$ and a Heisenberg superalgebra}
   \label{sec_heis}

We see from Theorem~\ref{th_main} that $\FG (X)$ has the
same size of the tensor product of the Fock space of an infinite-dimensional
Heisenberg algebra of rank $\dim K^0_G(X)$ and that of
an infinite-dimensional
Clifford algebra of rank $\dim K^1_G(X)$. It is our next step
to actually construct such an Heisenberg/Clifford algebra.
We will simply refer to as Heisenberg superalgebra from now on.

The dual of $\underline{K}_G(X)$, denoted by $\underline{K}_G(X)^*$,
is naturally $\Z_2$-graded as
identified with $\underline{K}^0_G(X)^* \bigoplus \underline{K}^1_G(X)^*$.
Denote by $\langle \cdot, \cdot \rangle$ the pairing
between $\underline{K}_G(X)^*$ and $\underline{K}_G(X)$.
For any $n, m \geq 1$ and $\eta \in \underline{K}_G(X)^*$,
we define an additive map
\begin{eqnarray}  \label{eq_ann}
  a_{-m} (\eta) : \underline{K}_{G_n}(X^n) \longrightarrow
            \underline{K}_{G_{n-m} }(X^{n-m} )
\end{eqnarray}
as the composition
\begin{eqnarray*}
  \underline{K}_{G_n}(X^n)
& \stackrel{\mbox{Res} }{\longrightarrow}
&              \underline{K}_{G_{m} \times G_{n -m} }(X^n)
\stackrel{k^{-1} }{\longrightarrow}
              \underline{K}_{G_{m}} (X^{m} )\bigotimes
               \underline{K}_{G_{n -m} }(X^{n -m} )     \\
&\stackrel{{ch}_m \otimes 1 }{\longrightarrow}
&    \underline{K}_{G} (X ) \bigotimes \underline{K}_{G_{n -m}} (X^{n -m} )
\stackrel{\eta \otimes 1 }{\longrightarrow}
              \underline{K}_{G_{n -m}} (X^{n -m} ).
\end{eqnarray*}

On the other hand, we define for any
$m \geq 1$ and $V \in \underline{K}_G (X)$ an additive map
\begin{eqnarray}  \label{eq_creat}
  a_{m} (V) : \underline{K}_{G_{n-m} }(X^{n-m} )
          \longrightarrow \underline{K}_{G_n}(X^n)
\end{eqnarray}
as the composition
\begin{eqnarray*}
  \underline{K}_{G_{n-m} }(X^{n-m} )
  &\stackrel{\omega_m (V) \boxtimes \cdot}{\longrightarrow} &
  \underline{K}_{G_m }(X^m) \bigotimes \underline{K}_{G_{n -m}} (X^{n -m} ) \\
  & \stackrel{k }{\longrightarrow} &
    \underline{K}_{G_m \times G_{n -m}}(X^n )
\stackrel{\mbox{Ind}}{\longrightarrow}
          \underline{K}_{G_n}(X^n).
\end{eqnarray*}

Let $\cal H$ be the linear span of the
operators $a_{-m} (\eta), a_{m} (V), m \geq 1$,
$\eta \in \underline{K}_G(X)^*,$ $V \in \underline{K}_G(X)$.
Clearly $\cal H$ admits a natural $\Z_2$-gradation
induced from that on $\underline{K}_G(X)$ and $ \underline{K}_G(X)^* $.
Below we shall use $[ \cdot , \cdot ]$ to denote the
supercommutator as well. It is understood that $[a, b]$ is
the anti-commutator $ab +ba$ when $a, b \in \cal H$ are both odd elements
according to the $\Z_2$-gradation.

\begin{theorem}    \label{th_heisenberg}
 When acting on $\FG (X)$, $\cal H$ satisfies the
 Heisenberg superalgebra commutation relations,
namely for $m , l \geq 1,$ $ \eta, \eta^{'} \in \underline{K}_G(X)^*, $
$V,W \in \underline{K}_G(X)$,
\begin{eqnarray}
  [ a_{-m} (\eta), a_{l} (V)] & =& l \delta_{m,l} \langle \eta, V \rangle ,
           \label{eq_heis}  \\
  {[ a_{m} (W)    , a_{l} (V)]} & =& 0,      \label{eq_induct} \\
  {[ a_{-m} ({\eta}), a_{-l} ({\eta}^{'}) ] }& = & 0.  \label{eq_restrict}
\end{eqnarray}
Furthermore $\FG (X)$ is an irreducible representation of the
Heisenberg superalgebra.
\end{theorem}

\begin{demo}{Proof}
 We may assume that $V, \eta$ are homogeneous, say of
 degree $v$ and $e$ where $v, e \in \{0, 1\}$, according to the
 $\Z_2$-grading of $\underline{K}_G(X)$ and its dual.
 We keep using the notations in the proof of Theorem~\ref{th_hopf}.

 Given $E \in \underline{K}_{G_r} (X^r)$, we first observe by the definitions
 (\ref{eq_ann}) and (\ref{eq_creat}) that $a_{-m}(\eta) a_l (V) E$
 (resp. $(-1)^{ve} a_l (V) a_{-m}(\eta) E$) is given by the composition
 from the top to the bottom along the left (resp. right)
 side of the diagram below:
 \begin{eqnarray*}
  \begin{array}{rcl}
                 &     \underline{K}_{G_r} (X^r)       &  \\
                 &     \downarrow  \!{ \omega_l (V) \boxtimes \cdot}    &  \\
                & \underline{K}_{G_l} (X^l)
                  \bigotimes  \underline{K}_{G_r} (X^r) &  \\
  {\footnotesize  Ind}\! \swarrow      &                               &
                             \searrow\!{\footnotesize  1 \otimes \Res } \\
 \underline{K}_{G_N} (X^N)   & & \underline{K}_{G_l} (X^l)
         \bigotimes \underline{K}_{G_m} (X^m)
         \bigotimes \underline{K}_{G_{r -m}} (X^{r -m})   \\
   \Vert  \;\;\;\quad        & & \quad  \downarrow T^{(l,m)} \otimes 1 \\
    \underline{K}_{G_N} (X^N)   & & \underline{K}_{G_m} (X^m)
         \bigotimes \underline{K}_{G_l} (X^l)
         \bigotimes \underline{K}_{G_{r -m}} (X^{r -m})   \\
  {\footnotesize \Res}\!{\searrow}   &                              &
          {\footnotesize P}  %
  {\swarrow}\!{\footnotesize  \Ind \otimes 1 }  \\
       & \underline{K}_{G_m} (X^m) \bigotimes \underline{K}_{G_n} (X^n) &  \\
                 & \;\;\;\downarrow
                    \!{\footnotesize  ch_m  \otimes 1}    &  \\
       & \underline{K}_{G} (X)  \bigotimes \underline{K}_{G_n} (X^n)     &  \\
           & \downarrow  \!{\footnotesize \eta  \otimes 1}  &  \\
                 &  \underline{K}_{G_n} (X^n)     &
 \end{array}
 \end{eqnarray*}
 Here and below it is understood that when a negative integer
 appears in indices the corresponding term is $0$.
 To simplify notations, we put $\Res$ (resp. $\Ind$) instead of
 the composition $k^{-1} \circ \Res$ (resp. $\Ind \circ k$)
 in the above diagram and below.

  We denote by ${\cal M}'$ the set of all the
 $2 \times 2$ matrices of the form (\ref{matr})
 satisfying (\ref{eq_condition}) except the following two matrices
\begin{eqnarray*}
 \left[ \begin{array}{cc}
  0        & m    \\
  l      & r -m
 \end{array} \right],
 \quad
 \left[ \begin{array}{cc}
  m        & 0    \\
  l-m      & r
 \end{array} \right].
 \end{eqnarray*}
 As in the proof of Theorem~\ref{th_hopf}, we apply Lemma~\ref{lem_ser}
 to the case $Y = X^N, \Gamma = G_N, H = G_l \times G_r,$
 and $L = G_m \times G_n$, where $l + r = m +n = N$.
 \begin{eqnarray}
  && \Res_{(m, n)} \Ind_{(l, r)}^N (\omega_l (V) \boxtimes E ) \nonumber \\
  &= &
  \bigoplus_{A \in \cal M} \Ind^{(m,n)}_A
  (1_{a_{11} } \otimes T^{(a_{12}, a_{21})} \otimes 1_{a_{22} } )
  (\Res_{A'} (\omega_l (V) \boxtimes E ) )  \nonumber \\
  & = &
   \bigoplus_{A \in \cal M} \Ind^m_{(a_{11}, a_{12})} F_1 \boxtimes
   \Ind^n_{(a_{21}, a_{22})} F_2   \nonumber \\
 & = &
   \bigoplus_{A \in {\cal M}'} \Ind^m_{(a_{11}, a_{12})} F_1 \boxtimes
   \Ind^n_{(a_{21}, a_{22})} F_2   \nonumber \\
  & & \bigoplus F_1 \boxtimes \Ind^n_{(l, r-m)} F_2
     \bigoplus F_1 \boxtimes \Ind^n_{(l-m, r)} F_2   \nonumber \\
  & = &
   \bigoplus_{A \in {\cal M}'} \Ind^m_{(a_{11}, a_{12})} F_1 \boxtimes
   \Ind^n_{(a_{21}, a_{22})} F_2  \nonumber \\
   & & \bigoplus ( 1_m \otimes \Ind_{(l, r- m)}^n )
   ( T^{(l,m)} \otimes 1_{l -m} )
   ( 1_l \otimes \Res_{(m, r -m)}) (\omega_l (V) \boxtimes E )
            \label{eq_mess}    \\
   & &     \bigoplus (1_m \otimes \Ind_{(l -m, r)}^n) )
     (\Res_{(m, l -m)} \otimes 1_r ) (\omega_l (V) \boxtimes E ) . \nonumber
 \end{eqnarray}

 We get $0$ when applying the map $ch_m \otimes 1 $ to
 $\Ind^m_{(a_{11}, a_{12})} F_1 \boxtimes \Ind^n_{(a_{21}, a_{22})} F_2$
 for $A \in {\cal M}'$  in (\ref{eq_mess})
 by Lemma~\ref{lem_orth}. When applying
 $( \eta \otimes 1 ) \circ ( ch_m \otimes 1 )$ to the
 second term of the r.h.s. of (\ref{eq_mess}), we obtain
 $  ( -1)^{ve} a_l (V) a_{-m}(\eta) E$. When applying $ch_m \otimes 1$ to
 the third term of the r.h.s. of (\ref{eq_mess}),
 we get $0$ if $ m \neq l$ by Lemma~\ref{lem_orth}.
 In the case $m =l$, the third term of the r.h.s.
 of (\ref{eq_mess}) is simply $\omega_l (V) \boxtimes E$. When applying
 $(\eta \otimes 1 ) \circ ( ch_m \otimes 1)$ to it,
 we get $l \langle \eta, V \rangle E$
 by Proposition~\ref{prop_property}. Putting all these
 pieces together, we have proved Eq.~(\ref{eq_heis}).

 We may assume that $W$ is homogeneous of degree $w \in \{0, 1\}$
 according to the $\Z_2$-grading of $\underline{K}_G(X)$.
 Eq.~(\ref{eq_induct}) is a consequence
 of the transitivity of the induction functor:
 \begin{eqnarray*}
 &&\Ind_{G_m \times G_{l +r} }^{G_{m + l +r} }
  \left( \omega_m (W) \boxtimes \Ind_{G_l \times G_r}^{G_{l +r} }
         ( \omega_l (V) \boxtimes E)
  \right)    \\
  &= & \Ind_{G_m \times G_{l +r} }^{G_{m + l +r} }
  \left( \Ind_{G_m \times G_l \times G_r}^{G_m \times G_{l +r} }
         (\omega_m (W) \boxtimes  \omega_l (V) \boxtimes E)
  \right)    \\
  &= & \Ind_{G_m \times G_l \times G_r}^{G_{m + l +r} }
         (\omega_m (W) \boxtimes  \omega_l (V) \boxtimes E)    \\
  &= & ( -1)^{vw}\Ind_{G_l \times G_{m +r} }^{G_{m + l +r} }
  \left( \Ind_{G_l \times G_m \times G_r}^{G_l \times G_{m +r} }
          ( \omega_l (V) \boxtimes \omega_m (W) \boxtimes E)
  \right)  \\
  &= & \Ind_{G_l \times G_{m +r} }^{G_{m + l +r} }
  \left( \omega_l (V) \boxtimes \Ind_{G_m \times G_r}^{G_{m +r} }
          ( \omega_m (W) \boxtimes E)
  \right).
  \end{eqnarray*}

 Similarly Eq.~(\ref{eq_restrict}) is a consequence
 of the transitivity of the restriction functor.
 The irreducibility of ${\cal F}_G(X)$ as a representation
 of $\cal H$ follows immediately from the
 $q$-dimension formula for ${\cal F}_G(X)$ given in
 Theorem~\ref{th_main}.
\end{demo}

In the special case $G =1$, the Heisenberg superalgebra
was constructed by Segal \cite{S2} which differs slightly
from ours. Our proof follows his strategy of proof as well.

\begin{remark}  \rm
 One may consider the enlarged space
 $$
  V_G := \FG (pt) \bigotimes \C [R(G)]
 $$
 where $\C [R(G)]$ is the group
 algebra of the lattice $R(G)$. Note that $V_G$ is the underlying
 space for a lattice vertex algebra \cite{B, FLM}.
 In particular, when $G$ is a finite
 subgroup of $SL_2(\Bbb C)$,
 the space $V_G$ is closely related to the Frenkel-Kac-Segal
 vertex representation of an affine Lie algebra. In this way, we
 are able to obtain a new link between the subgroups of $SL_2(\Bbb C)$
 and the affine Lie algebras widely known as
 the McKay correspondence.
 Connections among symmetric functions,
 $\FG (pt)$, $V_G$, and vertex operators will be developed
 in a forthcoming paper.

 More generally, one may consider
\begin{eqnarray}  \label{eq_latt}
  V_G(X) : = {\cal F}_G(X) \bigotimes \C [ K_G(X) ]
 \end{eqnarray}
 when $K_G(X)$ is torsion-free, where $\C [ K_G(X) ] $
 is the group algebra of the lattice $\C [ K_G(X) ]$
 (if $K_G(X)$ is not torsion-free
 we replace $K_G(X)$ in (\ref{eq_latt}) by the free
 part of $K_G(X)$ over $\Z$).
\end{remark}
\section{The orbifold Euler characteristic $e(X^n, G_n)$}
    \label{sec_orbi}
  In the study of string theory on orbifolds,
Dixon {\em et al} \cite{DHVW}
came up with a notion of {\em orbifold Euler characteristic} defined
as follows:
$$e (X, G) = \frac1{|G|} \sum_{g_1 g_2 = g_2 g_1} e(X^{g_1, g_2}),$$
where $X$ is a smooth $G$-manifold.
$X^{g_1, g_2}$ denotes the common fixed point
set of $g_1$ and $g_2$ and $e(\cdot)$ denotes
the usual Euler characteristic.
One easily shows \cite{HH} that the orbifold Euler characteristic
can be equivalently defined as
\begin{eqnarray}   \label{eq_euler}
  e (X, G) = \sum_{ [g] \in G_*} e (X^g / Z_G(g) ) .
\end{eqnarray}

Denote by $X^{(n)}$ the $n$-th symmetric product of $X$.
Recall that Macdonald's formula \cite{M}
relates $e( X^{(n)})$ to $e(X)$ as follows:
\begin{eqnarray}
 \sum_{n =0} e( X^{(n)}) q^n = ( 1 -q)^{- e(X)}. \label{eq_mac}
\end{eqnarray}

The following theorem relates the orbifold Euler
characteristic $e(X^n, G_n)$ to $e(X, G)$.
\begin{theorem} \label{th_orbi} We have
   $\sum_{n \geq 0} e(X^n, G_n) q^n =  \prod_{ r =1}^{\infty}
  ( 1 -q^r)^{- e(X, G)} .$
\end{theorem}

\begin{demo}{Proof}
  For an alternative proof see Remark~\ref{rem_another} below.
 By the definition of the orbifold Euler characteristic,
 Lemmas~\ref{lem_cent} and Lemma~\ref{lem_fix}, we have
  \begin{eqnarray*}
    \sum_{n \geq 0} e(X^n, G_n) q^n
   &= & \sum_{ n \geq 0} \sum_{ [a] \in (G_n)_*}
         e \left( (X^n)^a / Z_{G_n} (a)
           \right) q^n  \quad\quad\quad\quad\quad\quad
       \mbox{ by Eq. }(\ref{eq_euler}),   \\
   &\stackrel{(A)}= & \sum_{ n \geq 0} \sum_{\sum_r r m_r(c) = n}
         \prod_{c \in G_*}
        e \left( (X^c /Z_G(c))^{m_r(c)}
          \right)  q^n     \label{eq_ele}   \\
   &= & \prod_{c \in G_*} \prod_{r \geq 1}
   \left(
        \sum_{ m_r(c) \geq 0}
        e \left( (X^c /Z_G(c))^{m_r(c)}
          \right)  q^{m_r(c)}
   \right)               \nonumber  \\
   &= &
         \prod_{c \in G_*} \prod_{r \geq 1} ( 1 -q^r)^{ - e(X^c /Z_G(c))}
        \;\;   \mbox{ by applying }(\ref{eq_mac}) \mbox{ to } X^c /Z_G(c),
         \label{eq_macd} \\
   &= & \frac1{ \prod_{ r =0}^{\infty}
  ( 1 -q^r)^{  e(X, G)} }
    \;\;\quad\quad\quad\quad\quad\quad\quad\quad\quad\quad\quad
    \mbox{ by Eq. }(\ref{eq_euler}).
 \end{eqnarray*}
 Here Eq.~(A) follows from
 Lemma~\ref{lem_cent} and Lemma~\ref{lem_fix}.
\end{demo}

 In the case when $G$ is trivial, we recover a formula given in \cite{HH}.
\begin{remark}   \rm  \label{rem_another}
 According to Atiyah and Segal \cite{AS} the orbifold Euler
 characteristic can be calculated in terms of
 equivariant K-theory:
 \begin{eqnarray}   \label{eq_diff}
   e(X, G)= \dim K^0_G (X) - \dim K^1_G (X).
 \end{eqnarray}
  Theorem~\ref{th_orbi} follows from Theorem~\ref{th_main}
 by applying Eq.~(\ref{eq_diff}) to $K_{G_n} (X^n)$.
\end{remark}

One is interested \cite{DHVW, HH} in finding a resolution of
singularities
$$\widetilde{X/G} \longrightarrow X/G$$
with the property
$$
e(X, G) = e(\widetilde{X/G}).
$$
We assume that $X$ is a smooth
quasi-projective surface with such a property in the
following discussions. Denote by $X^{[n]}$ the Hilbert
scheme of $n$ points on $X$.
According to G\"ottsche \cite{G}, the Euler characteristic
of $X^{[n]}$ is given by:
\begin{equation} \label{eq_got}
 \sum_{n \geq 0} e( X^{[n]}) q^n
 = \prod_{ r =0}^{\infty}  ( 1 -q^r)^{- e(X)}.
\end{equation}

 We note that $X^n/ G_n $ is naturally identified with $ (X/G)^n /S_n$.
The following commutative diagram
 \begin{eqnarray*}
  \begin{array}{ccc}
       {\widetilde{X/G} }^{[n]}  & \rightarrow
                                      & \left(\widetilde{X/G}\right)^n  /S_n \\
               \downarrow        &             & \downarrow  \\
           X^n/ G_n              &  \equiv     & (X/G)^n /S_n
  \end{array}
 \end{eqnarray*}
implies that the Hilbert scheme
${\widetilde{X/G} }^{[n]}$ is a resolution of singularity of $X^n/ G_n$
(indeed it is semismall).
By applying Eq.~(\ref{eq_got}) to $\widetilde{X/G}$ and
comparing with Theorem~\ref{th_orbi}, we have the following corollary.
\begin{corollary}  \label{cor_gen}
  Let $X$ be a smooth quasi-projective surface and assume that
 there exists a smooth resolution ${\widetilde{X/G} }$ of singularities
 of the orbifold $X/G$ such that $e(\widetilde{X/G} ) = e(X, G)$.
 Then there exists a resolution of singularities of $X^n/ G_n$
 given by $ {\widetilde{X/G} }^{[n]}$ satisfying
 \begin{eqnarray*}
   e(X^n, G_n) = e \left( {\widetilde{X/G} }^{[n]} \right).
 \end{eqnarray*}
\end{corollary}

  The assumption of the existence of the
resolution ${\widetilde{X/G} }$ of singularities
of $X/G$ above is necessary as this is the special case of
$X^n /G_n$ for $n =1$.
In the setting of Corollary~\ref{cor_gen}, we conjecture
that ${\widetilde{X/G} }^{[n]} $ is a {\em crepant}
resolution of $X^n /G_n$ provided that
${\widetilde{X/G} }$ is a crepant resolution of singularities
of $X/G$.

We consider a special case in detail.
Let $X$ be the complex plane $\C^2$
acted upon by a finite subgroup $G$ of $SL_2(\Bbb C)$.
Via the McKay correspondence \cite{Mc}, there
is a one-to-one correspondence between the finite subgroups of
$SL_2(\Bbb C)$ and the Dynkin diagrams of simply-laced
types $A_n$, $D_n$, $E_6,$ $E_7$ and $E_8$. Let us denote
by $\frak g$ the simple Lie algebra corresponding to $G$.
From the exact correspondence we know that
the rank of $\frak g$ is $|G_*| -1$.

The quotient $\C^2/ G$ has an isolated simple singularity at $0$. There
exists a minimal resolution $\widetilde{\C^2/ G }$ of $\C^2/ G$
well known as ALE-spaces (cf. e.g. \cite{N}).
It is known that the second homology group $H_2 (\widetilde{\C^2/ G }, \Z)$
is isomorphic to the root lattice of $\frak g$, cf. e.g. \cite{N}.
In particular
\begin{equation} \label{eq_wonder}
  \dim  K \left(\widetilde{\C^2/ G } \right)
   = \dim  H \left(\widetilde{\C^2/ G } \right) = |G_*|.
\end{equation}
So if we apply Theorem~\ref{th_main} to the case $\widetilde{\C^2/ G }$
with a trivial group action, we have
\begin{eqnarray*}
  \sum_{n \geq 0} \dim K_{S_n}
  \left( \widetilde{\C^2/ G }^{n} \right) q^n
  =  \prod_{r \geq 1}  (1 - q^r)^{ - |G_*|} .
\end{eqnarray*}

On the other hand,
by the Thom isomorphism \cite{S1} we have
$$
 K_{G_n} (X^n) \approx K_{G_n} (pt ) = R(G_n).
$$
It follows from
the part~3) of Remark~\ref{rem_spec} and Eq. (\ref{eq_wonder}) that
\begin{eqnarray*}
  \sum_{n \geq 0} \dim {K}_{G_n} \left( \C^{2n} \right) q^n
  =  \prod_{r \geq 1}  (1 - q^r)^{ - |G_*|} .
\end{eqnarray*}

We can obtain another numerical coincidence from a
somewhat different point of view as follows.
For a general quasi-projective smooth surface $Y$, a well-known result
of Fogarty says that the Hilbert
scheme of $n$ points on $Y$, denoted by $Y^{[n]}$,
is a smooth $2n$ dimensional manifold.
The Betti numbers of  $Y^{[n]}$ were computed by G\"ottsche \cite{G}.
In particular, G\"ottsche's formula yields the dimension of the
homology group of $Y^{[n]}$:
\begin{eqnarray} \label{eq_dim}
 \sum_{n \geq 0} \dim H( Y^{[n]}) q^n
 = \prod_{r \geq 1}  (1 - q^r)^{ - \dim H(Y)}.
\end{eqnarray}
We apply Eq. (\ref{eq_dim}) to the case $Y = \widetilde{\C^2/ G }$.
It follows from Eq. (\ref{eq_wonder}) that
\begin{eqnarray*}
 \sum_{n \geq 0} \dim H( \widetilde{\C^2/ G }^{[n]}) q^n
 = \prod_{r \geq 1}  (1 - q^r)^{ - |G_*|}.
\end{eqnarray*}

Therefore we have proved the following.
\begin{proposition}  \label{prop_iden}
 The spaces $ {K}_{G_n} \left( \C^{2n} \right) \bigotimes \C$,
 $ {K}_{S_n} ( \widetilde{\C^2/ G }^{n} )  \bigotimes \C$
 and $ H( \widetilde{\C^2/ G }^{[n]}) $ have the same dimension.
\end{proposition}

Since the minimal resolution $\widetilde{\C^2 /G}$ of $\C^2 /G$
has no odd dimensional homology, we have (cf. \cite{HH, N})
$$
  e(\widetilde{\C^2 /G}) = \dim H(\widetilde{\C^2 /G})
  = e (\C^2, G) =|G_*|.
$$
The following corollary
is a special and important case of Corollary~\ref{cor_gen}.
\begin{corollary}    \label{cor_elem}
 Let $G$ be a finite subgroup of $SL_2(\Bbb C)$
 and $\widetilde{\C^2 /G}$ be the minimal resolution of $\C^2 /G$.
 Then the Hilbert scheme of $n$ points on $\widetilde{\C^2 /G}$
 is a resolution of singularities of $\C^{2n} /G_n$ whose Euler
 characteristic is equal to the orbifold Euler
 characteristic $e(\C^{2n}, G_n)$.
\end{corollary}

\begin{remark} \rm  \label{rem_dim}
 The fact that the (graded) dimension of $K_{S_n} \left( X^{n} \right)$
 equals that of the homology group of the Hilbert scheme
 of $n$ points of $X$ for a more general surface $X$ holds
 by the same argument as above. Bezrukavnikov and Ginzburg
 \cite{BG} have proposed a way to establish
 a direct isomorphism between these two spaces for an
 algebraic surface $X$.
 M.~de Cataldo and L.~Migliorini has recently established an isomorphism
 independently for any complex surface \cite{CM}.

 Proposition~\ref{prop_iden} can be generalized as follows.
 Assume that $X$ is a quasi-projective surface acted upon by $G$ and
 there exists a smooth resolution of singularities $\widetilde{X/ G }$
 of $X/G$ such that the dimension of $K_G^i (X)$ $(i =0, 1)$
 equals that of the even (resp. odd) dimensional homology
 group of $\widetilde{X/ G }$.
 Then we conclude that the dimension
 of $K_{G_n} \left( X^n \right)$ is equal to that
 of the homology group of the Hilbert scheme of $n$ points
 on $\widetilde{X/ G }$.
 We conjecture the existence of a natural isomorphism from
 $ {K}_{G_n} \left( X^{n} \right) \bigotimes \C$
 to $ H( \widetilde{ X/ G }^{[n]}) $, assuming the (necessary)
 existence of an isomorphism from $K_G (X) \bigotimes \C$
 to $ H( \widetilde{ X/ G }) $ or $K( \widetilde{ X/ G }) \bigotimes \C$.

 We believe that this is just a first
 indication of intriguing relations between the Hilbert
 scheme of $n$ points on $\widetilde{X/ G }$
 and the wreath product $G_n$.
 We will elaborate on this in another occasion.
 When $G$ is trivial, it reduces to
 well-known relations between the Hilbert scheme
 of $n$ points and the symmetric group $S_n$ (cf.~\cite{N, BG}).
\end{remark}

\begin{remark} \rm
 Let us consider a special case of Corollary~\ref{cor_elem}
 by putting $G = \Z_2$. The wreath product $G_n$ in this case is exactly
 the Weyl group of $B_n$ or $C_n$.
 It is interesting to compare with a ``Hilbert scheme''
 associated to a reductive group of type $B_n$ (or $C_n$)
 defined in \cite{BG}.
\end{remark}

\noindent{\bf Acknowledgement.} The starting point of
this work is the insight of Segal \cite{S2}.
I am grateful to him for his permission of using his unpublished results.
I am also grateful to Igor Frenkel who
first emphasized to me the importance of \cite{S2}
and to Mikio Furuta whose generous help is indispensable
for me to understand \cite{S2} for inspiring and fruitful discussions.
I thank Mark de Cataldo for helpful conversations
on Hilbert schemes. I also thank the referee for
his comments and suggestions which help to
improve the presentation of the paper.
This paper is a modified version of my MPI preprint
``Equivariant K-theory and wreath products'' in August 1998.
It is a pleasure to acknowledge
the warm hospitality of the Max-Planck Institut f\"ur Mathematik at Bonn.

\vspace{.1in}

\noindent{\bf Note 1 added.} The idea here of relating the
representation rings of wreath products associated to finite
groups of $SL_2 ( \C)$ and vertex representations of affine Lie
algebras has been fully developed in a paper ``Vertex
representations via finite groups and the McKay correspondence''
by I.~Frenkel, N.~Jing and the author. \vspace{.1in}

\noindent{\bf Note 2 added.} The connections among Hilbert
schemes, wreath products and K-theory have been developed in my
recent paper "Hilbert schemes, wreath products, and the McKay
correspondence".

\noindent
Department of Mathematics, Yale University, New Haven, CT 06520, USA \\
E-mail: wqwang@math.yale.edu

\end{document}